\documentclass[10pt]{amsart}

\usepackage[english]{babel}
\usepackage[T1]{fontenc}
\usepackage[utf8]{inputenc}
\usepackage{amsmath}
\usepackage{amssymb}
\usepackage{amsthm}
\usepackage{color}
\usepackage[a4paper]{geometry}
\usepackage{tikz}
\usepackage{mathdots}
\usepackage[all]{xy}
\usepackage{hyperref}
\usepackage{accents}
\usepackage{lmodern}
\usepackage{graphicx}

\theoremstyle{definition}
\newtheorem{definition}{Definition}[section]

\theoremstyle{plain}
\newtheorem{theorem}[definition]{Theorem}
\newtheorem{prop}[definition]{Proposition}

\newcommand{\mb}{\mathbb}
\newcommand{\mc}{\mathcal}

\newcommand{\bs}{\boldsymbol}
\newcommand{\Nm}{\textup{Nm}}
\newcommand{\Vol}{\textup{Vol}}

\begin{document}

\title{Sums of reciprocals of fractional parts}

\author{Reynold Fregoli}
\address{Department of Mathematics\\ 
Royal Holloway, University of London\\ 
TW20 0EX Egham\\ 
UK}
\email{Reynold.Fregoli.2017@live.rhul.ac.uk}

\subjclass{Primary ; Secondary }
\date{\today, and in revised form ....}

\dedicatory{}


\begin{abstract}
Let $\bs{\alpha}\in \mb{R}^N$ and $Q\geq 1$. We consider  the sum $\sum_{\bs{q}\in [-Q,Q]^N\cap\mb{Z}^N\backslash\{\bs{0}\}}\|\bs{\alpha}\cdot\bs{q}\|^{-1}$. 
Sharp upper bounds are known when $N=1$, using continued fractions or the three distance theorem. However, these techniques do not seem to apply in higher dimension.
We introduce a different approach, based on a general counting result of Widmer for weakly admissible lattices, to establish sharp upper bounds for arbitrary $N$.
Our result also sheds light on a question raised by L\^{e} and Vaaler in 2013 on the sharpness of their lower bound $\gg Q^N\log Q$. 
\end{abstract}

\maketitle

\section{Introduction}

Sums of reciprocals of fractional parts have been studied by many authors in the light of their tight connections to, e.g., 
uniform distribution theory, metric Diophantine approximation, and lattice point counting (see \cite{Beck:ProbDioph}, \cite{HarLi:SomeProblems1}, \cite{HarLi:SomeProblems2}, \cite{Kruse:Estimates}, \cite{KuipersNiederr:UnifDistr}, \cite{Schmidt:MetricalThms} and \cite{Shapira:ASolution}). In this note, we establish a new, sharp upper bound for such sums. 

We denote by $\|x\|$ the distance to the nearest integer of a real $x$. We write $A\ll B$ to mean that there exists a constant $c>0$ (absolute or depending only on the parameters indicated) such that $A\leq cB$. We use $|\cdot|_{2}$ to denote the Euclidean norm on $\mb{R}^{n}$ and $|\cdot|_{\infty}$ to denote the maximum norm.

Let $N\in \mb{N}:=\{1,2,3,\ldots\}$, and let $\bs{x}\cdot\bs{y}:=\sum_{i=1}^N x_i y_i$ be the standard inner product on $\mb{R}^N$.
Let $\bs{\alpha}$ be in $\mb{R}^{N}$ with $1,\alpha_1,\ldots,\alpha_N$ linearly independent over $\mb{Q}$, and suppose $Q_{1},\dotsc,Q_{N}\in(0,+\infty)$. L\^{e} and Vaaler \cite[Corollary 1.2]{LeVaaler:Sumsof} proved that for $X:=[-Q_1,Q_1]\times\cdots\times[-Q_N,Q_N]$ and $Q:=(Q_1\cdots Q_N)^{1/N}\geq 1$ it holds
 \begin{alignat}1\label{LV2}
Q^N\log Q\ll \sum_{\bs{q}\in X\cap\mb{Z}^N\backslash\{\bs{0}\}}\|\bs{\alpha}\cdot \bs{q}\|^{-1}.
\end{alignat}
They also showed that, whenever $\bs{\alpha}$ is multiplicatively badly approximable (see \cite{Bugeaud:Multiplicative} and \cite[(2.5)]{LeVaaler:Sumsof}), inequality (\ref{LV2}) is sharp\footnote{ Note that L\^{e} and Vaaler assumed $X:=[0,Q_1]\times\cdots\times[0,Q_N]$, but this makes no difference. Suppose indeed 
that the sum taken over, e.g., the set $[-Q_1,0]\times[0,Q_2]\times\cdots\times[0,Q_N]$ is big, then, we can multiply the first coordinate of $\bs{\alpha}$ by $-1$, obtaining a multiplicatively badly approximable vector such that the sum over $X$ is now big.}, i.e.,
\begin{alignat}1\label{LV3}
\sum_{\bs{q}\in X\cap\mb{Z}^N\backslash\{\bs{0}\}}\|\bs{\alpha}\cdot \bs{q}\|^{-1}\ll_{\bs{\alpha}}Q^N\log Q
\end{alignat}
for all $Q\geq 2$. However, if $N\geq 2$, every such $\bs{\alpha}$ yields a counterexample to the Littlewood conjecture.
This follows from two facts:
\begin{itemize}
\item a matrix is multiplicatively badly approximable if and only if its transpose is multiplicatively badly approximable (see \cite[Corollary 1]{German:Transference} and \cite[Theorem 2.2]{LeVaaler:Sumsof});
\item every submatrix of a multiplicatively badly approximable matrix is itself multiplicatively badly approximable (see \cite[(2.6)]{LeVaaler:Sumsof}).
\end{itemize}
To date, Littlewood's conjecture is still open and the set of counterexamples to it, if not empty, has been shown to be extremely sparse. Indeed, Einsiedler, Katok, and Lindenstrauss \cite{EinseidlerKatokLindenstrauss:Invariant} have proved that its Hausdorff dimension is zero. Despite this, it can be shown that for some vectors in $\mb{R}^{N}$ we can reverse inequality (\ref{LV2}). 

\begin{definition}
\label{def:badlyapp}
Let $N\in\mb{N}$ and $\bs{\alpha}\in\mb{R}^{N}$. Let $\phi:(0,+\infty)\to(0,1]$ be a non-increasing\footnote{ We say that $\phi$ is non-increasing if $\phi(x)\geq\phi(y)$, whenever $x<y$.}, real-valued function. We say that $\bs{\alpha}$ is a $\phi$-badly approximable vector if
$$|\bs{q}|_{\infty}^{N}\left\|\bs{\alpha}\cdot\bs{q}\right\|\geq\phi(|\bs{q}|_{\infty})$$
for all $\bs{q}\in\mb{Z}^{N}\setminus\{\bs{0}\}$. If $\phi$ can be chosen constant, we simply say that $\bs{\alpha}$ is badly approximable.
\end{definition}

It was kindly pointed out to me by Victor Beresnevich that a simple gap-principle, already known in the literature (see \cite[Proof of Lemma 3.3, p.123]{KuipersNiederr:UnifDistr}), proves (\ref{LV3}) for all badly approximable vectors $\bs{\alpha}$ in the special case $X=[-Q,Q]^N$. Note that the set of such vectors has full Hausdorff dimension in $\mb{R}^{N}$ (see \cite{Schmidt:BadlyApprox}). We recall this proof here. Let $Q\geq 2$ and suppose that $\bs{\alpha}\in\mb{R}^{N}$ is $\phi$-badly approximable. Then, for all distinct $\bs{q}_{1}, \bs{q}_{2}\in\mb{Z}^{N}\cap X$, we have  
\begin{alignat}1\label{www}
\|\bs{\alpha}\cdot\bs{q}_{1}\pm\bs{\alpha}\cdot\bs{q}_{2}\|=\|\bs{\alpha}\cdot(\bs{q}_{1}\pm\bs{q}_{2})\|\geq \frac{\phi(2Q)}{(2Q)^{N}}.
\end{alignat}
It follows that 
$$\left|\|\bs{\alpha}\cdot\bs{q}_{1}\|-\|\bs{\alpha}\cdot\bs{q}_{2}\|\right|\geq \frac{\phi(2Q)}{(2Q)^{N}}.$$
Therefore, none of the intervals 
$$\left[0,\frac{\phi(2Q)}{(2Q)^{N}}\right),\left[\frac{\phi(2Q)}{(2Q)^{N}},\frac{2\phi(2Q)}{(2Q)^{N}}\right),\left[\frac{2\phi(2Q)}{(2Q)^{N}},\frac{3\phi(2Q)}{(2Q)^{N}}\right),\ldots$$ 
contains more than one number of the form $\|\bs{\alpha}\cdot \bs{q}\|$ ($\bs{q}\in \mb{Z}^{N}\cap X$), and no such number lies in the first interval. Hence,
\begin{alignat}1\label{????}
\sum_{\bs{q}\in X\cap\mb{Z}^N\backslash\{\bs{0}\}}\|\bs{\alpha}\cdot \bs{q}\|^{-1}\leq \sum_{j=1}^{|X\cap\mb{Z}^N\backslash\{\bs{0}\}|}\frac{(2Q)^{N}}{j\phi(2Q)}\ll_N \frac{Q^N\log Q}{\phi(2Q)}.
\end{alignat}
For $N=1$,  the above upper bound can be improved on to 
\begin{alignat}1\label{??}
\ll Q\log Q+\frac{Q}{\phi(Q)}
\end{alignat}
using the theory of continued fractions (see \cite[Theorem 2, p.37-40]{Lang:IntrotoDiophApprox}). 
The bound  given in (\ref{??}) is best possible. However,   
if one allows the upper bound to be expressed in terms of the (least) denominator  $q_{K}$ of the $K$-th convergent of $\alpha$, 
an even more precise result can be obtained via the three distance theorem, as shown by Beresnevich and Leong \cite[Corollary 1]{BeresnevichLeong:Somsof}.

Nevertheless, neither the techniques based on continued fractions nor those of Beresnevich and Leong using the three distance theorem seem to generalise in an obvious way
to higher dimension.
In this note we introduce yet another method, based on a recent counting result for weakly admissible lattices, which allows us to extend (\ref{??}) to arbitrary dimension $N$.

\begin{theorem}
\label{thm:theorem3}
Let $X:=[-Q,Q]^N$ and let $\bs{\alpha}\in\mb{R}^{N}$ be a $\phi$-badly approximable vector. Then, we have
\begin{equation}
\label{eq:theorem3-1}
Q^{N}\log\left(Q\phi(Q)\right)\ll_{N}\sum_{\substack{\mathbf{q}\in X \cap\mb{Z}^{N}\setminus\{\bs{0}\}}}\left\|\bs{\alpha}\cdot\bs{q}\right\|^{-1}\ll_{N}Q^{N}\log Q+\frac{Q^{N}}{\phi(Q)}
\end{equation}
for all $Q\geq 1$.
\end{theorem}

Clearly, the lower bound is superseded by (\ref{LV2}) (and non trivial only if $\log(Q\phi(Q))\geq 1$), but we have decided to include it in Theorem (\ref{thm:theorem3}) nonetheless. This is because both bounds in (\ref{eq:theorem3-1}) can be proved almost at once and the method we use is substantially different from L\^{e} and Vaaler's.
Note that the upper bound, in conjunction with the Khintchine-Groshev Theorem, implies that for every $\varepsilon>0$ the set of $\bs{\alpha}\in \mb{R}^N$ for which
$$\sum_{\substack{\mathbf{q}\in X\cap \mb{Z}^{N}\setminus\{\bs{0}\}}}\left\|\bs{\alpha}\cdot\bs{q}\right\|^{-1}\ll_{\bs{\alpha}}Q^{N}(\log Q)^{1+\varepsilon},$$
has full Lebesgue measure. Moreover, the upper bound in (\ref{eq:theorem3-1}) is sharp, in the sense that for all $N\geq 1$ there exists a sequence of positive integers $Q_{i}\to +\infty$ such that
\begin{equation}
\label{eq:sharpness}
\sum_{\substack{\mathbf{q}\in X_{i}\cap \mb{Z}^{N}\setminus\{\bs{0}\}}}\left\|\bs{\alpha}\cdot\bs{q}\right\|^{-1}\gg Q_{i}^{N}\log Q_{i}+\frac{Q_{i}^{N}}{\phi(Q_{i})},
\end{equation}
where $X_{i}:=[-Q_{i},Q_{i}]^{N}$. To show this, let $\phi$ be chosen maximal, i.e.,
$$\phi(x)=\min\{|\bs{q}|_{\infty}^{N}\|\bs{\alpha}\cdot\bs{q}\|\ |\ 0<|\bs{q}|_{\infty}\leq x,\ \bs{q}\in\mb{Z}^{N}\}$$  
for $x\geq 1$ and $\phi(x)=1$ for $x<1$. Let $\bs{q}_{i}\in\mb{Z}^{N}$ be a sequence of pairwise distinct vectors such that $\phi(|\bs{q}_{i}|_{\infty})=|\bs{q}_{i}|_{\infty}^{N}\|\bs{\alpha}\cdot\bs{q}_{i}\|$, and set $Q_{i}:=|\bs{q}_{i}|_{\infty}$. Then, $\|\bs{\alpha}\cdot\bs{q}_{i}\|^{-1}=Q^{N}_{i}/\phi(Q_{i})$. This implies that
\begin{equation}
\label{eq:sharpness2}
\sum_{\substack{\mathbf{q}\in X_{i}\cap \mb{Z}^{N}\setminus\{\bs{0}\}}}\left\|\bs{\alpha}\cdot\bs{q}\right\|^{-1}\geq\frac{Q_{i}^{N}}{\phi(Q_{i})}.
\end{equation}
Hence, (\ref{eq:sharpness}) follows from (\ref{LV2}) and (\ref{eq:sharpness2}).

The main tool to prove Theorem \ref{thm:theorem3} is Proposition \ref{prop:mainestimate1}, which gives a precise estimate for the size of the set
\begin{alignat}1
\label{eq:M1}
M(\bs{\alpha},\varepsilon, Q):=\left\{(p,\bs{q})\in\mb{Z}^{1+N}\setminus\{\bs{0}\}\ \right|\ \left.\left|\bs{\alpha}\cdot\bs{q}+p\right|\leq\varepsilon,\ |\bs{q}|_{\infty}\leq Q\right\},
\end{alignat}
where $0<\varepsilon\leq 1/2$ and $Q\geq 1$.
\begin{prop}
\label{prop:mainestimate1}
Let $\bs{\alpha}\in\mb{R}^{N}$ be a $\phi$-badly approximable vector. Let $Q\geq 1$ and $0<\varepsilon\leq 1/2$. Then, we have
$$\left|\left|M(\bs{\alpha},\varepsilon, Q)\right|-2^{N+1}\varepsilon Q^{N}\right|\ll_{N}\left(\frac{\varepsilon Q^{N}}{\phi(Q)}\right)^{\frac{N}{N+1}}.$$
\end{prop}
With Proposition \ref{prop:mainestimate1} at hand, Theorem \ref{thm:theorem3} can be derived from a simple dyadic summation. Moreover,  Proposition \ref{prop:mainestimate1} is a straightforward consequence of a recent, general lattice point counting
result, due to Widmer \cite[Theorem 2.1]{Widmer:WeakAdmiss}, which we recall in the next section.

\section{Weakly admissible lattices and counting lattice points}

To state \cite[Theorem 2.1]{Widmer:WeakAdmiss}, we need to introduce some notation. We follow the notation used in \cite{Widmer:WeakAdmiss}, except that we write $\mc{N}$ (instead of $N$) for $\sum_{i=1}^{n}m_{i}$.
We assume throughout that $\mc{N}\geq 2$.

Let $n$ be a positive integer and let $\mc{S}:=(\bs{m},\bs{\beta})$, where $\bs{m}:=(m_{1},\dotsc,m_{n})\in\mb{N}^{n}$ and $\bs{\beta}:=(\beta_{1},\dotsc,\beta_{n})\in(0,+\infty)^{n}$. Let
$$\Nm_{\bs{\beta}}\left(\underline{\bs{x}}\right):=\prod_{i=1}^{n}\left|\bs{x}_{i}\right|_{2}^{\beta_{i}}$$
be the multiplicative norm induced by $\bs{\beta}$ on $\mb{R}^{m_{1}}\times\dotsb\times\mb{R}^{m_{n}}\ni\underline{\bs{x}}:=(\bs{x}_{1},\dotsc,\bs{x}_{n})$. Let  $I$ be a non-empty subset of $\{1,\dotsc,n\}$ and let also
$$C=C_{I}:=\left\{\left.\underline{\bs{x}}\in\mb{R}^{m_{1}}\times\dotsb\times\mb{R}^{m_{n}}\ \right|\ \bs{x}_{i}=\bs{0}\ \ \forall i\in I\right\}.$$
We fix the couple $(\mc{S},C)$ and for any $\Gamma\subset\mb{R}^{m_{1}}\times\dotsb\times\mb{R}^{m_{n}}$ we consider the quantity
$$\nu\left(\Gamma,\varrho\right):=\inf\left\{\left. \Nm_{\bs{\beta}}(\underline{\bs{x}})^{\frac{1}{t}}\ \right|\ \underline{\bs{x}}\in\Gamma\setminus C,\ \left|\underline{\bs{x}}\right|_{2}<\varrho\right\},$$
where $t:=\beta_{1}+\dotsb+\beta_{n}$. We observe that $\nu(\Gamma,\cdot)$ is a decreasing function of $\varrho$, bounded by below from $0$. Hence, we have the following definition from \cite{Widmer:WeakAdmiss}.
\begin{definition}
\label{def:weakadmiss}
A full rank lattice $\Lambda$ in $\mb{R}^{m_{1}}\times\dotsb\times\mb{R}^{m_{n}}$ is said to be weakly admissible for the couple $(\mc{S},C)$, if $\nu\left(\Lambda,\varrho\right)>0$ for all $\varrho\in(0,+\infty)$.
\end{definition}

\noindent Before stating the counting theorem, we require some more notation. For any $\Gamma\subset\mb{R}^{m_{1}}\times\dotsb\times\mb{R}^{m_{n}}$ we define
$$\lambda_{1}\left(\Gamma\right):=\inf\left\{\ |\bs{x}|_{2}\ |\ \bs{x}\in\Gamma\setminus\{\bs{0}\}\right\},$$
and we set
$$\mu\left(\Gamma,\varrho\right):=\min\left\{\lambda_{1}\left(\Gamma\cap C\right),\nu\left(\Gamma,\varrho\right)\right\}.$$

We can now state a simplified version of \cite[Theorem 2.1]{Widmer:WeakAdmiss}, streamlined for our application.
\begin{theorem}[Widmer]
\label{thm:2.1}
Let $n\in\mb{N}$ and let $(\mc{S},C)$ be a couple as in Definition (\ref{def:weakadmiss}). For \\
$\bs{Q}:=\left(Q_{1},\dotsc,Q_{n}\right)\in(0,+\infty)^{n}$ we set
$$\overline{Q}:=\left(Q_{1}^{\beta_{1}}\dotsm Q_{n}^{\beta_{n}}\right)^{\frac{1}{t}}$$
and $Q_{\max}:=\max\{Q_{1},\dotsc,Q_{n}\}$. Let
$$Z_{\bs{Q}}:=\prod_{i=1}^{n}[-Q_{i},Q_{i}]^{m_{i}}\subset\mb{R}^{m_{1}}\times\dotsb\times\mb{R}^{m_{n}}$$
and let $\Lambda$ be a weakly admissible lattice for the couple $(\mc{S},C)$. Then, there exists a real constant $c=c(\mc{N})>0$, only depending on the quantity $\mc{N}:=\sum_{i=1}^{n}m_{i}$, such that
$$\left|\left|Z_{\bs{Q}}\cap\Lambda\right|-\frac{\Vol\left(Z_{\bs{Q}}\right)}{\det\Lambda}\right|\leq c\!\!\!\!\!\!\displaystyle\inf_{\ \ \ \ 0<B\leq Q_{\max}}\left(\frac{\overline{Q}}{\mu(\Lambda,B)}+\frac{Q_{\max}}{B}\right)^{\mc{N}-1},$$
where $\Vol(Z_{\bs{Q}})$ denotes the volume of the set $Z_{\bs{Q}}$ and $\det\Lambda$ denotes the determinant of the lattice $\Lambda$.
\end{theorem}

\section{Proofs}

\subsection{Proof of Proposition \ref{prop:mainestimate1}}
First we note that if $(p,\bs{q})\in M(\bs{\alpha},\varepsilon,Q)$, then
\begin{equation}
\label{eq:emptycase}
\varepsilon Q^{N}\geq\left\|\bs{q}\cdot\bs{\alpha}\right\||\bs{q}|_{\infty}^{N}\geq\phi(Q).
\end{equation}
Suppose that $\varepsilon Q^{N}/\phi(Q)< 1$. Then, $M(\bs{\alpha},\varepsilon,Q)=\emptyset$, by (\ref{eq:emptycase}). Moreover,
$$2^{N+1}\varepsilon Q^{N}\ll_{N}\varepsilon Q^{N}\leq\frac{\varepsilon Q^{N}}{\phi(Q)}\leq\left(\frac{\varepsilon Q^{N}}{\phi(Q)}\right)^{\frac{N}{N+1}}.$$
Hence, Proposition \ref{prop:mainestimate1} holds true whenever $\varepsilon Q^{N}/\phi(Q)< 1$. From now on, we can assume that
\begin{equation}
\label{eq:notemptycase}
\frac{\varepsilon Q^{N}}{\phi(Q)}\geq 1.
\end{equation}
Let
$$ A_{\bs{\alpha}}:=\left(\begin{array}{@{}c|ccc@{}}
    1 & \alpha_{1} & \dots & \alpha_{N} \\\hline
    0 &  &  & \\
    \vdots & & \phantom{\scriptstyle{N}}\scalebox{1.5}{$\text{I}$}_{N} & \\
    0 &  &  & 
  \end{array}\right).$$
Define $\Lambda_{\bs{\alpha}}:=A_{\bs{\alpha}}\mb{Z}^{N+1}\subset\mb{R}^{N+1}$ and $Z_{\varepsilon,Q}:=[-\varepsilon,\varepsilon]\times[-Q,Q]^{N}\subset\mb{R}^{N+1}$. Then,
\begin{equation}
\label{eq:cardinalities}
\left|M(\bs{\alpha},\varepsilon, Q)\right|=\left|\Lambda_{\bs{\alpha}}\cap Z_{\varepsilon,Q}\right|-1,
\end{equation}
since $\bs{0}\in\Lambda_{\bs{\alpha}}\cap Z_{\varepsilon,Q}$.  Therefore, to prove Proposition \ref{prop:mainestimate1}, it suffices to estimate the quantity $\left|\Lambda_{\bs{\alpha}}\cap Z_{\varepsilon,Q}\right|$. To this end, we use Theorem \ref{thm:2.1}.

Let $n=2$ and let $\bs{m}=\bs{\beta}:=(1,N)$. Let $C:=C_{I}$, with $I:=\{2\}$. Then, all vectors $\bs{v}\in\Lambda\setminus C$ have the form
$$\bs{v}=A_{\bs{\alpha}}\begin{pmatrix} p \\ \bs{q} \end{pmatrix}=\begin{pmatrix} \bs{\alpha}\cdot\bs{q}+p \\ \bs{q} \end{pmatrix},$$
where $\bs{q}\in\mb{Z}^{N}\setminus\{\bs{0}\}$ and $p\in\mb{Z}$.
Recall now that $\bs{\alpha}$ is a $\phi$-badly approximable vector. Hence, for all $\bs{v}\in\Lambda_{\bs{\alpha}}\setminus C$ it holds
$$\Nm_{\bs{\beta}}(\bs{v})=|\bs{\alpha}\cdot\bs{q}+p||\bs{q}|_{2}^{N}\geq\|\bs{\alpha}\cdot\bs{q}\||\bs{q}|_{\infty}^{N}\geq\phi\left(\left|\bs{q}\right|_{\infty}\right)\geq\phi(|\bs{v}|_{2}).$$
Since $\phi$ is non-increasing, we can conclude that
\begin{equation}
\label{eq:weakadmiss}
\nu\left(\Lambda_{\bs{\alpha}},\varrho\right)\geq\phi(\varrho)^{\frac{1}{N+1}}>0
\end{equation}
for all real $\varrho\geq\lambda_{1}(\Lambda_{\bs{\alpha}}\setminus C)$. However, if $\varrho<\lambda_{1}\left(\Lambda_{\bs{\alpha}}\setminus C\right)$, then $\nu(\Lambda_{\bs{\alpha}},\varrho)= +\infty$ and (\ref{eq:weakadmiss}) trivially holds true. This shows that $\Lambda_{\bs{\alpha}}\subset\mb{R}\times\mb{R}^{N}$ is weakly admissible for the couple $((\bs{m},\bs{\beta}),C)$. We can thus apply Theorem \ref{thm:2.1}, with $\Lambda=\Lambda_{\bs{\alpha}}$ and $Z_{\bs{Q}}=Z_{\varepsilon,Q}$. By choosing $B:=Q_{\max}=Q$, we get
\begin{equation}
\label{eq:estimate}
\left|\left|Z_{\varepsilon,Q}\cap\Lambda_{\bs{\alpha}}\right|-\frac{\Vol\left(Z_{\varepsilon,Q}\right)}{\det\Lambda_{\bs{\alpha}}}\right|\leq c\left(\frac{\left(\varepsilon Q^{N}\right)^{\frac{1}{N+1}}}{\mu(\Lambda,Q)}+1\right)^{N}.
\end{equation}
Since $\det\Lambda_{\bs{\alpha}}=1$ and $\Vol\left(Z_{\varepsilon,Q}\right)=2^{N+1}\varepsilon Q^{N}$, to conclude the proof, we just need to estimate the right-hand side of (\ref{eq:estimate}). We observe that $\lambda_{1}(\Lambda_{\bs{\alpha}}\cap C)=\lambda_{1}\left(\mb{Z}\times\{\bs{0}\}\right)=1$. Hence,
$$\mu(\Lambda,Q)\geq\min\left\{1,\phi(Q)^{\frac{1}{N+1}}\right\}=\phi(Q)^{\frac{1}{N+1}},$$
by (\ref{eq:weakadmiss}). Combining (\ref{eq:cardinalities}) and (\ref{eq:estimate}), and using (\ref{eq:notemptycase}), we find
$$\left|\left|M(\bs{\alpha},\varepsilon, Q)\right|-2^{N+1}\varepsilon Q^{N}\right|\ll_{N}\left(\left(\frac{\varepsilon Q^{N}}{\phi(Q)}\right)^{\frac{1}{N+1}}+1\right)^{N}\ll_{N}\left(\frac{\varepsilon Q^{N}}{\phi(Q)}\right)^{\frac{N}{N+1}}.$$
This completes the proof. 

\subsection{Proof of Theorem \ref{thm:theorem3}}

We start by observing that
\begin{align}
\label{eq:coreq1}
\sum_{\substack{\mathbf{q}\in[-Q,Q]^{N}\\ \cap\ \mb{Z}^{N}\setminus\{\bs{0}\}}}\left\|\bs{\alpha}\cdot\bs{q}\right\|^{-1} & \leq\sum_{k=1}^{\infty}2^{k+1}\left|\left\{\left.(p,\bs{q})\in\mb{Z}^{N+1}\setminus\{\bs{0}\}\ \right| \right.\right.\ 2^{-k-1}<\left|\bs{\alpha}\cdot\bs{q}+p\right|
\left.\left.\leq 2^{-k},\ \left|\bs{q}\right|_{\infty}\leq Q\right\}\right|\nonumber \\
& \leq\sum_{k=1}^{\infty}2^{k+1}\left|M(\bs{\alpha},2^{-k}, Q)\right|=\sum_{k=1}^{\left\lfloor\log_{2}(Q^{N}/\phi(Q))\right\rfloor}2^{k+1}\left|M(\bs{\alpha},2^{-k}, Q)\right|,
\end{align}
where the last equation is due to (\ref{eq:emptycase}).
Now, by Proposition \ref{prop:mainestimate1}, we know that
\begin{equation}
\label{eq:asymptotics}
\left|\left|M(\bs{\alpha},2^{-k}, Q)\right|-2^{N+1-k}Q^{N}\right|\ll_{N}\left(\frac{2^{-k}Q^{N}}{\phi(Q)}\right)^{\frac{N}{N+1}}.
\end{equation}
Hence, (\ref{eq:coreq1}) yields
\begin{align}
\sum_{\substack{\mathbf{q}\in[-Q,Q]^{N}\\ \cap\ \mb{Z}^{N}\setminus\{\bs{0}\}}}\left\|\bs{\alpha}\cdot\bs{q}\right\|^{-1} & \ll_{N}\sum_{k=1}^{\left\lfloor\log_{2}(Q^{N}/\phi(Q))\right\rfloor}2^{k+1}\left(2^{N+1-k}Q^{N}+\left(\frac{2^{-k}Q^{N}}{\phi(Q)}\right)^{\frac{N}{N+1}}\right)\nonumber \\
 & \ll_{N}\sum_{k=1}^{\left\lfloor\log_{2}\left(Q^{N}/\phi(Q)\right)\right\rfloor}\left(Q^{N}+2^{\frac{k}{N+1}}\left(\frac{Q^{N}}{\phi(Q)}\right)^{\frac{N}{N+1}}\right)\nonumber \\[10pt] 
 & \ll_{N}Q^{N}\log_{2}\left(\frac{Q^{N}}{\phi(Q)}\right)+\left(\frac{Q^{N}}{\phi(Q)}\right)^{\frac{1}{N+1}}\left(\frac{Q^{N}}{\phi(Q)}\right)^{\frac{N}{N+1}}\label{eq:coreq1.5} \\[15pt] 
 & \ll_{N} Q^{N}\log Q+\frac{Q^{N}}{\phi(Q)}\label{eq:coreq1.75},
\end{align}
where (\ref{eq:coreq1.5}) follows from the trivial estimate $\sum_{k=1}^{K}2^{k/(N+1)}\leq 2^{K/(N+1)+1}$ and (\ref{eq:coreq1.75}) is due to the fact that $1/\phi(Q)\geq\log\left(1/\phi(Q)\right)$. This proves the upper bound.

To prove the lower bound, we notice that
\begin{align}
\label{eq:coreq2}
\sum_{\substack{\mathbf{q}\in[-Q,Q]^{N}\\ \cap\ \mb{Z}^{N}\setminus\{\bs{0}\}}}\left\|\bs{\alpha}\cdot\bs{q}\right\|^{-1} & \geq\sum_{k=1}^{\infty}2^{k}\left|\left\{\left.(p,\bs{q})\in\mb{Z}^{N+1}\setminus\{\bs{0}\}\ \right| \right.\right.\ 2^{-k-1}<\left|\bs{\alpha}\cdot\bs{q}+p\right|
\left.\left.\leq 2^{-k},\ \left|\bs{q}\right|_{\infty}\leq Q\right\}\right|\nonumber \\
 & \geq\sum_{k=1}^{\infty}2^{k}\left(\left|M(\bs{\alpha},2^{-k}, Q)\right|-\left|M(\bs{\alpha},2^{-k-1}, Q)\right|\right).
\end{align}
From Proposition \ref{prop:mainestimate1}, we also know that for all $k\geq 1$ and $Q\geq 1$
\begin{equation}
\label{eq:coreq3}
\left|\left|M(\bs{\alpha},2^{-k}, Q)\right|-2^{N+1-k}Q^{N}\right|\leq c_{N}\left(\frac{2^{-k}Q^{N}}{\phi(Q)}\right)^{\frac{N}{N+1}},
\end{equation}
where $c_{N}$ is a positive constant. Hence, whenever $1\leq k\leq\log_{2}\left(Q^{N}\phi(Q)^{N}/c_{N}^{N+1}\right)=:K$, we have
$$\left(2^{N+1}-1\right)2^{-k}Q^{N}\leq \left|M(\bs{\alpha},2^{-k},Q)\right|\leq\left(2^{N+1}+1\right)2^{-k}Q^{N}.$$
This, in turn, shows that
\begin{equation}
\label{eq:coreq5}
\left|M(\bs{\alpha},2^{-k},Q)\right|-\left|M(\bs{\alpha},2^{-k-1},Q)\right|\geq 2^{-k}Q^{N},
\end{equation}
when $1\leq k\leq K-1$. Therefore, provided $K-1\geq 1$, we can plug (\ref{eq:coreq5}) into (\ref{eq:coreq2}) and restrict the sum to $k\leq\left\lfloor K-1\right\rfloor$. This yields the lower bound
\begin{equation}
\label{eq:coreq6.5}
\sum_{\substack{\mathbf{q}\in[-Q,Q]^{N}\\ \cap\ \mb{Z}^{N}\setminus\{\bs{0}\}}}\left\|\bs{\alpha}\cdot\bs{q}\right\|^{-1}\geq\lfloor K-1\rfloor Q^{N}\geq (K-2)Q^{N},
\end{equation} 
which, of course, remains true for $K-1< 1$. Now, a trivial lower bound for the sum of the reciprocals is
$$\sum_{\substack{\mathbf{q}\in[-Q,Q]^{N}\\ \cap\ \mb{Z}^{N}\setminus\{\bs{0}\}}}\left\|\bs{\alpha}\cdot\bs{q}\right\|^{-1}\gg_{N}Q^{N},$$
since $\left\|\bs{\alpha}\cdot\bs{q}\right\|^{-1}\geq 2$ for all $\bs{q}\neq 0$. Hence, recalling that $K=\log_{2}\left(Q^{N}\phi(Q)^{N}/c_{N}^{N+1}\right)$, we conclude that
$$\sum_{\substack{\mathbf{q}\in[-Q,Q]^{N}\\ \cap\ \mb{Z}^{N}\setminus\{\bs{0}\}}}\left\|\bs{\alpha}\cdot\bs{q}\right\|^{-1}\gg_{N}\left(K-2\right)Q^{N}+\log_{2}\left(4c_{N}^{N+1}\right)Q^{N}\gg_N Q^{N}\log\left(Q\phi(Q)\right).$$

\section*{Acknowledgements}
I'm most grateful to my supervisor, Martin Widmer, for his encouragement and precious advice. 
It has been a pleasure to discuss some aspects of this note with Victor Beresnevich, who provided some very useful feedback on it. I would like to thank Royal Holloway, University of London, for funding my position here, and my office mates for their support and cheerful presence.

\addcontentsline{toc}{section}{\bibname}
\bibliographystyle{plain}
\bibliography{Bibliography}

\end{document}